\documentclass[12pt]{article}

\usepackage{amsmath,amsthm,amssymb}
\usepackage{amsfonts}

\newtheorem{Tm}{\bf  Theorem}

\newtheorem{La}{\bf  Lemma }
\newtheorem{Rk}{\bf  Remark }

\begin{document}

\title{The index of singular integral operators \\ 
with discontinuous  oscillating coefficients}
\author{ E.V. Akulich, A.V. Lebedev\\
\\
Belarus State University, \\
Skariny av. 4, 220035, Minsk, Belarus\\
Institute of Mathematics, Bia{\l}ystok University,\\ ul. Akademicka 2, 
PL-15-424, Bia{\l}ystok, Poland} \maketitle

\begin{abstract}
The article presents the procedure of the index calculation for the elements of the 
algebra generated by one dimensional singular integral operators with discontinuous 
oscillating coefficients. 
\end{abstract}

{\small {\ KEY WORDS: singular integral operator, 
oscillating coefficients, $C^\ast$-algebra,
Fredholm operator, index, local representative }

\medskip

{\bf 2000 Mathematics Subject Classification:} 45E05, 32A55


\ \\

In the work by I.B. Simonenko and I.Ts. Gohberg and N.Ja. Krupnik 
[\ref{Simonenko1}--\ref{GohbergKrupnik}] the Fredholm theory of
the singular integral operators with {\em piecewise continuous} coefficients was  
developed.  
This theory was based on the  {\em local representatives} calculus for the 
operators considered. For the singular integral operators with discontinuous 
{\em oscillating} 
 coefficients the local representatives calculus was obtained in  [\ref{N1}].
 In the present article we continue to investigate these objects and give a 
 procedure of the index calculation for  the operators under consideration.
 This procedure goes well for  a finite set of  points of oscillation 
 and in order not to overload the presentation we describe it for the situation 
 when the coefficients have the only one point of oscillation.

Let  $M$ be an oriented closed simple Lyapunov curve,
 $B$ be the $C^*$-algebra generated by operators acting on 
 $L^2(M)$ and having the form 
$$
b=c_1+c_2S,
$$
where  $S$ is the singular integral operator and  $c_1$ and  $c_2$ are
the operators of multiplication by functions from the space 
$C(M\setminus
\{m_0\})$ having finite limits for   $m\to m_0\pm 0$, $m_0\in
M$.

The primitive ideal space $\mathrm{Prim}\ B$ of the algebra 
 $B$ can be represented as the disjoint union 
$$
\mathrm{Prim}\ B=M_{+}\cup M_{-}\cup {\mathbf
{R}}_{m_0}\cup(m_0,+,\pm)\cup(m_0,-,\pm),
$$
where  $M_{+}$ and $M_{-}$ are the copies of the set $M\setminus \{m_0 \}$,
${\mathbf {R}}_{m_0}$ is the straight line, $(m_0,+,\pm)$, $(m_0,-,\pm)$
are four points. 
By  $(m,\pm)$ we denote the points of  $M_{\pm}$ respectively.
The base of the topology on  $\mathrm{Prim} \ B$ is defined in the following way:
[\ref{Plamenevski}],\ [\ref{IPart}, 33.12],
[\ref{IIPart2}, 53.12]

\begin{enumerate}
\item[(a)]  a neighbourhoud of  $(m_0,+,+)$ is the union of  sets 
$$
\biggl [ \bigcup\limits_{m\in [m_0,m_1)} (m,+)\biggr ]\cup \{t\in
\mathbf{R}_{m_0} : t>N\},
$$
where  $N\in \mathbf{R}$ and  $m_1$ is an arbitrary point such that 
$m_0<m_1$ (here  $<$ is the order defined by the orientation of 
 $M$);

\item[(b)] a neighbourhood of  $(m_0,-,+)$ is the union of sets 
$$
\biggl [ \bigcup\limits_{m\in (m_2,m_0]} (m,+)\biggr ]\cup \{t\in
\mathbf{R}_{m_0} : t<N\};
$$

\item[(c)] a neighbourhood of  $(m_0,-,-)$ is the union of sets
$$
\biggl [ \bigcup\limits_{m\in [m_0,m_1)} (m,-)\biggr ]\cup \{t\in
\mathbf{R}_{m_0} : t>N\};
$$

\item[(d)] a neighbourhood of  $(m_0,+,-)$ is the union of sets 
$$
\biggl [ \bigcup\limits_{m\in (m_2,m_0]} (m,-)\biggr ]\cup \{t\in
\mathbf{R}_{m_0} : t<N\};
$$

\item[(e)] a neighbourhood of   $t\in \mathbf{R}_{m_0}$ is 
an open interval on  $\mathbf{R}_{m_0}$ containing  $t$;

\item[(f)] a neighbourhood of   $(m,\pm)$ for  $m\neq m_0$ 
is defined in the standard way  (as a neighbourhood of a point on a curve).
\end{enumerate}

As it has been already mentioned the local representatives calculus and the 
Fredholm theory of the operators from the algebra $B$
was developed  in  [\ref{Simonenko1} --- \ref{GohbergKrupnik}].

The symbol of an element  $b\in B$ is defined in terms of its local representatives 
that are given in the following way: 
$$
b^{\pm}(m)=c_1(m)\pm c_2(m),\ m\in M_{\pm},
$$
$$
b(m_0,+,\pm)=\mathop {\lim }\limits_{m\to m_0+0}
b^{\pm}(m)=c_1(m_0+0)\pm c_2(m_0+0),
$$
$$
b(m_0,-,\pm)=\mathop {\lim }\limits_{m\to m_0-0}
b^{\pm}(m)=c_1(m_0-0)\pm c_2(m_0-0),
$$
$$
b(m_0)(t)=\left(%
\begin{array}{cc}
  {\varphi}_{11}(t) & {\varphi}_{12}(t) \\
  {\varphi}_{21}(t) & {\varphi}_{22}(t) \\
\end{array}%
\right),
$$
where  ${\varphi}_{11}(t),\ {\varphi}_{22}(t)\in C(\overline {\mathbf
R}),\ {\varphi}_{12}(t),\ {\varphi}_{21}(t)\in \overset{\mathrm
o}{C} (\overline  {\mathbf R})$ and  $C(\overline {\mathbf R})$
is the algebra of continuous functions having limits when    $t\to \pm
\infty$ and  $\overset{\mathrm o}{C}(\overline  {\mathbf R})$ 
is the algebra of continuous functions tending to zero when  $t\to \pm
\infty$.

Moreover:
$$
b(m_0)(+\infty)=\begin{pmatrix}
  b(m_0,+,+) & 0 \\
  0 & b(m_0,-,-)
\end{pmatrix},\ \ \
$$
$$
b(m_0)(-\infty)=\begin{pmatrix}
  b(m_0,+,-) & 0 \\
  0 & b(m_0,-,+)
\end{pmatrix}.
$$

It is known [\ref{Simonenko1}--\ref{GohbergKrupnik}]
  that  an operator  $b\in B$ is Fredholm iff 
\begin{itemize}
  \item[(1)] $b^{\pm}(m)\ne 0 \ \forall m\in M_{\pm}$;
  \item[(2)] $\det \left[ b(m_0)(t)\right]\ne 0 \ \forall
  t\in {\mathbf{R}}_{m_0}$;
  \item[(3)] \text{ there exist nonzero limits }
$$
\mathop{\lim}\limits_{t \to\pm\infty}\det
\left[b(m_0)(t)\right]=\det\left[b(m_0)(\pm\infty)\right].
$$
\end{itemize}
And if (1)--(3) are satisfied then  ([\ref{GohbKrup1}], [\ref{GohbKrup2}])
\begin{equation}
{\rm {ind}}\ b= -\frac{1}{2\pi} \left( \arg
b^{+}(m){\arrowvert}_{M_{+}}- \arg b^{-}(m){\arrowvert}_{M_-}
+\arg\left[\det b(m_0)(t)\right] \arrowvert ^{+\infty}_{-\infty}
\right),\label{indB}
\end{equation}
where  $\left[\ \ \right]_{M\setminus\{m_0\}}$ is the increase 
of the function over the curve 
 $M\setminus\{m_0\}$ and  $\left[\ \
\right]|^{+\infty}_{t=-\infty}$ is the increase of the function over the straight line 
$\mathbf R$.

Let $h\in \mathbf R$ and  $U_h$ be the operator of multiplication by a function 
 $a_h(m)$ that  is continuous on 
 $M\setminus \{ m_0\}$ and such that on a certain symmetric neighbourhood 
 $O(m_0)$ of  $m_0\in M$ \,  $a_h(m)$  has the form 
\begin{equation}
\label{e1}
a_h(m)=\left\{%
\begin{array}{ll}
    e^{-\mathrm i h\ln (m_0-m)} & \hbox{where $m<m_0$,} \\
    e^{-\mathrm i h\ln (m-m_0)} & \hbox{where $m_0<m$,} \\
\end{array}%
\right.
\end{equation}
and $a_h(m)$ is equal to  1 out of $O(m_0)$.

We denote by  $C^{\ast}(B,U_h)$ the  $C^{\ast}$-algebra generated by the algebra 
 $B$ and the operators  $U_h, \ h\in \mathbf R$. The local representatives of elements of this 
 algebra were constructed in [\ref{N1}] and in terms of the local representatives 
 the conditions for the elements of $C^{\ast}(B,U_h)$
 to be Fredholm operators 
 were written out. 
 Let us write out the local representatives for the operator  $U_h$:
$$
(U_h)^{\pm}(m)=a_h(m),\ m\in M_{\pm},
$$
$$
U_h(m_0,+,\pm)f(t)=U_h(m_0,-,\pm)f(t)=U_h(m_0)f(t)=T_hf(t),
$$
where  $T_h$ is the shift operator  $T_hf(t)=f(t+h), \ \ f\in L^2 (\mathbf R , \mathbf C ^2)$. In this case the 
matrix function  $b(m_0)(t), t \in \mathbf R$ written out above is identified with the operator of multiplication by this matrix function in the space $L^2 (\mathbf R , \mathbf C ^2)$.

The aim of this article is to calculate the index  of  a  Fredholm operator 
 $d\in C^{\ast} (B,U_h)$ that has the form 
\begin{equation}\label{e2}
  d=b_0+b_1U_h,
\end{equation}
where  $b_i\in B,\ i=0,1,\ h\in \mathbf{R}$ in terms of the elements of the algebra 
 $B$, that is to give  a procedure of finding an operator  $d'\in B$,
such that 
$$
{\rm{ind\ }}d={\rm{ind\ }}d',
$$
so that one can apply formula  (\ref{indB}).

One of the principal technical instruments that are used 
in the procedure of the index calculation is a 
modification of theorem  46.20
[\ref{IIPart1}], which can be written 
out in terms of the objects of the present paper in the following way.  

\begin{Tm}
\label{t1}
\it{The operator  $d(m_0)=b_0(m_0)+b_1(m_0)T_h$ 
(that is the local representative of the operator  $d$  (\ref{e2})
on  ${\mathbf{R}}_{m_0}$)
 in the space  $L^2(\mathbf{R},{\mathbf{C}}^2)$ 
 is invertible iff there exist non degenerate continuous matrix functions 
  $w_1$ and $s_1$ (that are simultaneously diagonal at 
$\pm\infty$ or skew diagonal at one of the infinities and 
diagonal at the other infinity) such that 

\begin{equation}\label{e3}
    w^{-1}_1d(m_0)s_1=e_0+e_1T_h,
\end{equation}
where
$$
e_0=
\begin{pmatrix}
  I_l & 0 \\
  0 & e^0_{22}
\end{pmatrix},\
e_1=
\begin{pmatrix}
  e^1_{11} & 0 \\
  0 & I_{2-l}
\end{pmatrix},\ l=1,2,
$$
$$
r(e^0_{22}T^{-1}_h)<1,\ \ r(e^1_{11}T_h)<1,
$$
or 
\begin{equation}\label{e4}
    w^{-1}_1d(m_0)s_1=e_0+e_1T_h,
\end{equation}
where 
$$
e_0=
\begin{pmatrix}
  e^0_{11} & 0 \\
  0 & I_{2-l}
\end{pmatrix},\
e_1=
\begin{pmatrix}
  I_l & 0 \\
  0 & e^1_{22}
\end{pmatrix},\ l=1,2,
$$
$$
r(e^0_{11}T_h)<1,\ \ r(e^1_{22}T^{-1}_h)<1.
$$
(The operators  $e^0_{22}T^{-1}_h$ in (\ref{e3}) and  $e^1_{22}T^{-1}_h$ in 
(\ref{e4}) are considered on the invariant subspace defined by the projection 
$$
\begin{pmatrix}
  0 & 0 \\
  0 & I_{2-l}
\end{pmatrix},
$$
and the operators  $e^0_{11}T_h$ in (\ref{e3}) and  $e^1_{11}T_h$ in 
(\ref{e4}) are considered in the subspace defined by the projection 
$$
\begin{pmatrix}
  I_l & 0 \\
  0 & 0
\end{pmatrix},
$$
here $I_s$ is the identity operator acting on a space of dimension $s$.) 
Operator  $T_h$ is the shift operator:
$$
T_hf(t)=f(t+h).
$$ }\end{Tm}

The proof of this theorem can be obtained by the same scheme as the proof of 
theorem  46.20 [\ref{IIPart1}] by taking into account the explicit form of the operator 
 $d(m_0)$.

\

\begin{Rk}
\em The simultaneous diagonality 
(skew diagonality) of the matrices  $w_1$ and  $s_1$ at one and the same infinity 
follows from their construction 
(under the scheme of  the proof of theorem  46.20 [\ref{IIPart1}])
and from the fact that 
$b_i(m_0)(\pm\infty),\ i=0,1$ are the diagonal matrices.
\end{Rk}
\begin{Rk}
\em 
  We exclude from the considerations the situation when the matrices $w_1$ and $s_1$
are skew diagonal at both the infinities since in this case 
one can  take the diagonal at 
$\pm\infty$ matrices 
$$
w_2(t)=w_1(t)\left(%
\begin{array}{cc}
  0 & 1 \\
  1 & 0 \\
\end{array}%
\right)\ \ \hbox{and}\ \ \ s_2(t)=s_1(t)\left(%
\begin{array}{cc}
  0 & 1 \\
  1 & 0 \\
\end{array}%
\right).
$$
And in this way we obtain 
\begin{equation*}\begin{split}
w^{-1}_2d(m_0)s_2&=\left(%
\begin{array}{cc}
  0 & 1 \\
  1 & 0 \\
\end{array}%
\right)w^{-1}_1d(m_0)s_1\left(%
\begin{array}{cc}
  0 & 1 \\
  1 & 0 \\
\end{array}%
\right)=\\&=\left(%
\begin{array}{cc}
  0 & 1 \\
  1 & 0 \\
\end{array}%
\right)(e_0+e_1T_h)\left(%
\begin{array}{cc}
  0 & 1 \\
  1 & 0 \\
\end{array}%
\right)=e_0+e_1T_h.
\end{split}
\end{equation*}
\end{Rk}

\

A useful observation is the next
\begin{La} 
\label{l1}
\it {If the matrices  $w_1$ and $s_1$ are simultaneously diagonal at 
 $\pm\infty$ then there exist operators $s$ and  $w$ from  $B$ that have  zero indexes
and are such that for their local representatives
on ${\mathbf {R}}_{m_0}$ we have
$$
s(m_0)=s_1\ \ and \ \ w(m_0)=w_1 .
$$
}
\end{La}
{\bf Proof.} By conditions of the lemma it is clear that the matrix functions 
 $s_1$ and  $w_1$ can be continued by means of non vanishing continuous functions 
  $s^{\pm}$ and  $w^{\pm}$ on  $M_{\pm}$ in such a way that 
$$
\arg s^{+}(m){\arrowvert}_{M_{+}}- \arg s^{-}(m){\arrowvert}_{M_-}
=\arg\left[\det s_1(t)\right] \arrowvert ^{+\infty}_{-\infty},
$$
$$
\arg w^{+}(m){\arrowvert}_{M_{+}}- \arg w^{-}(m){\arrowvert}_{M_-}
=\arg\left[\det w_1(t)\right] \arrowvert ^{+\infty}_{-\infty}.
$$
After this is done it is enough to apply formula  (\ref{indB}). $\square $

\

We start our considerations with  the following  technical statement. 
\begin{La}
\label{l2} 
{\it Let $A(t), B(t), X(t), Y(t), Z_0(t), Z_1(t), \ t\in
\overline{\mathbf{R}}_{m_0}$ be  two dimensional continuous matrix functions 
such that $ X(t), Y(t)$  are non degenerate at  every 
 $t$ and $A(\pm\infty)$  and $B(\pm\infty)$ are diagonal.
 
Let also  $\infty$ denotes one of the 
(fixed) infinities: either   $+\infty$ or  $-\infty$ and 
$A(\infty)=[a_{ij}]^2_{i,j=1}$, $B(\infty)=[b_{ij}]^2_{i,j=1}$ 
and so on.

If in the space 
$L^2(\mathbf{R},{\mathbf{C}}^2)$ the operators of multiplication by matrix functions 
 $A, B, X, Y,
Z_0, Z_1$ satisfy the equality 
\begin{equation}\label{ABXYZ}
X(A+BT_h)Y=Z_0+Z_1T_h,
\end{equation}
 then:
\begin{itemize}
\item[(i)] the matrices  $Z_0(t)$ and  $Z_1(t)$ have the form 
$$
Z_0(t)=\left(%
\begin{array}{cc}
  z^0_{11}(t) & 0 \\
  0 & 1 \\
\end{array}%
\right),\ \ \
Z_1(t)=\left(%
\begin{array}{cc}
 1 & 0 \\
  0 & z^1_{22}(t) \\
\end{array}%
\right),
$$
where  $|z^0_{11}|<1$ and  $|z^1_{22}|<1$ iff 

$X(\infty)$ and  $Y(\infty)$  are simultaneously diagonal
and
\begin{equation}\label{i1}
|a_{11}|<|b_{11}|,\ |a_{22}|>|b_{22}|,\ a_{22}x_{22}y_{22}=1,\
b_{11}x_{11}y_{11}=1
\end{equation}
or 

$X(\infty)$ and $Y(\infty)$  are simultaneously skew diagonal 
and 
\begin{equation}\label{i2}
|a_{11}|>|b_{11}|,\ |a_{22}|<|b_{22}|,\ a_{11}x_{21}y_{12}=1,\
b_{22}x_{12}y_{21}=1;
\end{equation}

\item[(ii)] the matrices  $Z_0(t)$ and  $Z_1(t)$ have the form 
$$
Z_0(t)=\left(%
\begin{array}{cc}
  1 & 0 \\
  0 & z^0_{22}(t) \\
\end{array}%
\right),\ \ \
Z_1(t)=\left(%
\begin{array}{cc}
 z^1_{11}(t) & 0 \\
  0 & 1 \\
\end{array}%
\right),
$$
where  $|z^0_{22}|<1$ and  $|z^1_{11}|<1$ iff 

$X(\infty)$ and $Y(\infty)$  are simultaneously diagonal and
\begin{equation}
\label{ii1}
|a_{11}|>|b_{11}|,\ |a_{22}|<|b_{22}|,\ a_{11}x_{11}y_{11}=1,\
b_{22}x_{22}y_{22}=1
\end{equation}
or 

$X(\infty)$ and  $Y(\infty)$  are simultaneously skew diagonal and 
\begin{equation}
\label{ii2}
|a_{11}|<|b_{11}|,\ |a_{22}|>|b_{22}|,\ a_{22}x_{12}y_{21}=1,\
b_{11}x_{21}y_{12}=1;
\end{equation}

\item[(iii)] if 
$$
Z_0(t)=\left(%
\begin{array}{cc}
  1 & 0 \\
  0 & 1 \\
\end{array}%
\right),\ \ \
Z_1=\left(%
\begin{array}{cc}
  z^1_{11}(t) & 0 \\
  0 & z^1_{22}(t) \\
\end{array}%
\right),
$$
where  $|z^1_{11}|<1$ and  $|z^1_{22}|<1$ at least for one infinity 
then  $X(\pm\infty)$ and  $Y(\pm\infty)$ can be taken to be diagonal 
at both the infinities and 
\begin{equation}\label{iii}
|a_{11}|>|b_{11}|,\ |a_{22}|>|b_{22}|,\ a_{11}x_{11}y_{11}=1,\
a_{22}x_{22}y_{22}=1;
\end{equation}

\item[(iv)] if 
$$
Z_0(t)=\left(%
\begin{array}{cc}
  z^0_{11}(t) & 0 \\
  0 & z^0_{22}(t) \\
\end{array}%
\right),\ \ \
Z_1=\left(%
\begin{array}{cc}
  1 & 0 \\
  0 & 1 \\
\end{array}%
\right),
$$
where  $|z^0_{11}|<1$ and  $|z^0_{22}|<1$ at least at one infinity then 
$X(\pm\infty)$ and  $Y(\pm\infty)$ can be taken to be diagonal at both the 
infinities and 
\begin{equation}\label{iv}
|a_{11}|>|b_{11}|,\ |a_{22}|>|b_{22}|,\ b_{11}x_{11}y_{11}=1,\
b_{22}x_{22}y_{22}=1.
\end{equation}

\end{itemize}
}
\end{La}
{\bf Proof.}  Under the conditions of (i) the equality  (\ref{ABXYZ})
is equivalent to the system:
\begin{equation}
\label{Sistem i}
\left\{%
\begin{array}{ll}
    b_{11}x_{11}y_{11}+b_{22}x_{12}y_{21}=1,\\
    a_{11}x_{11}y_{12}+a_{22}x_{12}y_{22}=0,\\
    b_{11}x_{11}y_{12}+b_{22}x_{12}y_{22}=0,\\
    a_{11}x_{21}y_{11}+a_{22}x_{22}y_{21}=0,\\
    b_{11}x_{21}y_{11}+b_{22}x_{22}y_{21}=0,\\
    a_{11}x_{21}y_{12}+a_{22}x_{22}y_{22}=1,\\
    |a_{11}x_{11}y_{11}+a_{22}x_{12}y_{21}|<1,\\
    |b_{11}x_{21}y_{12}+b_{22}x_{22}y_{22}|<1.
\end{array}%
\right.
\end{equation}

By solving the second and the third equations of the system with respect to the unknowns 
$x_{11}y_{12}$ and  $x_{12}y_{22}$ we obtain that if 
$a_{11}b_{22}-a_{22}b_{11}\ne 0$ then 
$$
x_{11}y_{12}=x_{12}y_{22}=0.
$$
Under the same condition the forth and the fifth equations imply 
$$
x_{21}y_{11}=x_{22}y_{21}=0.
$$
If  $x_{11}=0$ then  $x_{12}\ne 0$, $x_{21}\ne 0$ (since the matrix 
$X(\infty)$ is non degenerate), so  $y_{11}=y_{22}=0$. 
Therefore  $y_{12}\ne 0$, $y_{21}\ne 0$ (since the matrix 
$Y(\infty)$ is non degenerate) and  $x_{22}=0$. Thus we have shown that the matrices 
 $X(\infty)$ and $Y(\infty)$ are skew diagonal. 

Substituting  $x_{11}=x_{22}=y_{11}=y_{22}=0$ into the first 
and the sixth equations and the inequalities of  system 
 (\ref{Sistem i}) we obtain 
(\ref{i2}).

If  $x_{11}\ne 0$ then one can show in an analogous way that  $X(\infty)$ 
and 
$Y(\infty)$ are diagonal and  (\ref{i1}) holds.

If  $\frac{a_{11}}{a_{22}}=\frac{b_{11}}{b_{22}}=\lambda\ne 0$,
then it follows from  (\ref{Sistem i}) that 
$$
\left\{%
\begin{array}{ll}
    |a_{22}||\lambda x_{11}y_{11}+x_{12}y_{21}|<1,  \\
    b_{22}(\lambda x_{11}y_{11}+x_{12}y_{21})=1,  \\
    |b_{22}||\lambda x_{21}y_{12}+x_{22}y_{22}|<1,  \\
     a_{22}(\lambda x_{21}y_{12}+x_{22}y_{22})=1.  \\
\end{array}%
\right.
$$
And we arrive at the inequalities 
$$
\big| \frac{a_{22}}{b_{22}}\big| <1 \ \ {\rm and }\ \  
\big| \frac{b_{22}}{a_{22}}\big| <1
$$
which is a contradiction.

The case  (ii) can be considered in an analogous way. 
The corresponding system has the form: 
$$
\left\{%
\begin{array}{ll}
    a_{11}x_{11}y_{11}+a_{22}x_{12}y_{21}=1,\\
    a_{11}x_{11}y_{12}+a_{22}x_{12}y_{22}=0,\\
    b_{11}x_{11}y_{12}+b_{22}x_{12}y_{22}=0,\\
    a_{11}x_{21}y_{11}+a_{22}x_{22}y_{21}=0,\\
    b_{11}x_{21}y_{11}+b_{22}x_{22}y_{21}=0,\\
    b_{11}x_{21}y_{12}+b_{22}x_{22}y_{22}=1,\\
    |b_{11}x_{11}y_{11}+b_{22}x_{12}y_{21}|<1,\\
    |a_{11}x_{21}y_{12}+a_{22}x_{22}y_{22}|<1.
\end{array}%
\right.
$$

In the case  (iii) we obtain the non degeneracy of the matrix 
 $A(t)$ for every 
$t\in \overline{\mathbf{R}}_{m_0}$ and the matrix 
$X(t)$ can be chosen to be equal to  
$A^{-1}(t)$ and  $Y(t)$ to be equal to the identity matrix. Thus 
(\ref{iii}) holds true.

The case  (iv) can be proved in the similar way. 
The proof is complete. $\square $

\

The properties of the topology of 
$\mathrm{Prim\ }B$ imply that not changing the index of a Fredholm 
 operator  $d$  of the form 
(\ref{e2}) one can change its coefficients and operators  $w$ and  $s$
in such a way that   $b_i(m)$, $w(m)$, $s(m)$ will be constant  (equal respectively to 
 $b_i(m_0\pm0)$, $w(m_0\pm 0)$, $s(m_0\pm 0)$) in a certain symmetric 
 neighbourhood  $O'(m_0)\subset O(m_0)$. In what follows we shall consider precisely such 
 'changed' operators.

When calculating the index of  $d$ we shall consider in sequel  all
the possible situations mentioned in theorem  \ref{t1}.

\begin{Tm}
\label{t2}
\it{A Fredholm operator $d$ of the form  (\ref{e2}) satisfies the inequalities \\
 (case I)
\begin{equation}
\label{caseI}
|b^{\pm}_0(m_0+0)|>|b^{\pm}_1(m_0+0)|,\
|b^{\pm}_0(m_0-0)|>|b^{\pm}_1(m_0-0)|
\end{equation}
or the inequalities  (case  II)
\begin{equation}
\label{caseII}
|b^{\pm}_0(m_0+0)|<|b^{\pm}_1(m_0+0)|,\
|b^{\pm}_0(m_0-0)|<|b^{\pm}_1(m_0-0)|)
\end{equation}
iff there exist  zero index operators  $w$ and  $s$ from 
the algebra  $B$ such that the local representative for 
$e=w^{-1}ds$ on the straight line  ${\mathbf {R}}_{m_0}$
is the operator  (\ref{e3}) (respectively  (\ref{e4})) for  $l=2$. 
In this case there exists a homotopy in the class of Fredholm operators 
of $C^\ast (B,U_h)$ of the  operator $e$ 
to an operator $d'\in B$ and therefore 
$$
\mathrm{ind}\ d =\mathrm{ind}\ e = \mathrm{ind}\ d'.
$$
}
\end{Tm}

{\bf Proof.} The proof exploits theorem \ref{t1} and we shall use the 
same notation.
The inequalities  (\ref{caseI}) and (\ref{caseII}) were obtained in lemma \ref{l2}  
for the cases  (iii) and  (iv) respectively  ($X(t)=w^{-1}_1(t),\
Y(t)=s_1(t),\ A(t)=b_0(m_0)(t),\ B(t)=b_1(m_0)(t),\ Z_0(t)=
\newline e_0(t),\ Z_1(t)=e_1(t)$).

By lemma \ref{l1} there exist zero index operators  $w$ and  $s$
such that 
 $w(m_0)(t)=w_1(t)$ and  $s(m_0)(t)=s_1(t)$.

In case  I it follows from (\ref{iii}) that 
$$
e^{\pm}(m_0+0)=1+\frac{b_1^{\pm}(m_0+0)}{b_0^{\pm}(m_0+0)}T_h,
$$
$$
e^{\pm}(m_0-0)=1+\frac{b_1^{\pm}(m_0-0)}{b_0^{\pm}(m_0-0)}T_h.
$$

Now we shall construct the homotopy of the operator $e$ to  
an operator $e' \in B$. 

To start with we define the homotopy of the operator 
$e(m_0)=e_0(m_0)+e_1(m_0)T_h$ (that is the homotopy of the corresponding
local representatives)
on  $\mathbf{R_{m_0}}\cup (m_0,+,\pm)\cup (m_0,-,\pm)$.

For every  $\tau\in[0,1]$ we set 
$$
(e'_{\tau})(m_0)(t)=\begin{pmatrix}
  1 & 0 \\
  0 & 1
\end{pmatrix}+\begin{pmatrix}
  (1-\tau)e^1_{11}(t) & 0 \\
  0 & (1-\tau)e^1_{22}(t)
\end{pmatrix}T_h,
$$$$
e'_{\tau}(m_0,+,\pm)=1+(1-\tau)
\frac{b_1^{\pm}(m_0+0)}{b_0^{\pm}(m_0+0)} T_h,
$$
$$
e'_{\tau}(m_0,-,\pm)=1+(1-\tau)
\frac{b_1^{\pm}(m_0-0)}{b_0^{\pm}(m_0-0)} T_h.
$$

For  $\tau=1$ we have 
$$
e'_1(m_0)(t)=e_0(m_0)(t)=\begin{pmatrix}
  1 & 0 \\
  0 & 1
\end{pmatrix},
$$
$$
e'_1(m_0,+,\pm)=e'_1(m_0,-,\pm)=1.
$$

Now let us extend this homotopy onto  $M_{+}$  and   $M_{-}$:
$$
(e'_{\tau})^{\pm}(m)=(w^{-1})^{\pm}(m)(b^{\pm}_0(m)+{\varphi}_{\tau}(m)b^{\pm}_1(m)a_h(m))s^{\pm}(m),
$$
where
$$
{\varphi}_{\tau}(m)=\left\{%
\begin{array}{ll}
    1-\tau\cdot\frac{m-m_1}{m_0-m_1} & \hbox{when  $m\in (m_0,m_1)$,} \\
    1-\tau\cdot\frac{m-m_2}{m_0-m_2} & \hbox{where $m\in (m_2,m_0)$,} \\
    1 & \hbox{where $m\notin (m_2,m_0)\cup(m_0,m_1)$,} \\
\end{array}%
\right.
$$
and $m_1$ is a certain point from the neighbourhood   on  $M_{+}$ and  $M_{-}$neighbourhood   on  $M_{+}$ and  $M_{-}$
corresponding to the neighbourhood  $O'(m_0)$ 
on  $M\setminus\{m_0\}$ and 
$m_2$  is the point that  is symmetric to  $m_1$ with respect to  $m_0$.

The operators  $(e'_{\tau})^{\pm}(m)$ are invertible for every 
$\tau\in[0,1]$ since in the opposite case 
in the neighbourhood corresponding to  $O'(m_0)$
one obtains 
$$
b^{\pm}_0(m)+{\varphi}_{\tau}(m)b^{\pm}_1(m)a_h(m)=0.
$$
And it follows that 
$$
|b^{\pm}_0(m)|=|{\varphi}_{\tau}(m)||b^{\pm}_1(m)|.
$$ 
But since the functions 
$b^{\pm}_i(m),\ i=1,2,$ are constants on this neighbourhood and 
$|{\varphi}_{\tau}(m)|\le 1$ this leads to a contradiction with 
(\ref{caseI}).

For every $\tau\in[0,1]$ the invertible operators  $(e'_{\tau})^{\pm}(m)$,
$e'_{\tau}(m_0,+,\pm)$, $e'_{\tau}(m_0,-,\pm)$ and  $e'_{\tau}(m_0)$
define the symbol of a certain Fredholm operator  $e'_{\tau}\in
C^{\ast}(B,U_h)$ (they are its local representatives).  Thus  
$$
\mathrm{ind}\ d=\mathrm{ind}\ e=\mathrm{ind}\ e'_0=\mathrm{ind}\
e'_1,
$$
but  $e'_1$ is an element of the algebra  $B$. Set $d'= e'_1$. The index of the latter operator can be 
calculated  by means of  (\ref{indB}).

In order to reduce case II to case I it is enough to multiply the 
operator  $d$ by  $U^{-1}_h$ (under this operation index does not change).
The theorem is proved. $\square $

\

To prove the next theorem we need the operator 
$P_1+Q_1U^{-1}_h$ where  $P_1$ and $Q_1$
are the operators of multiplication by the functions $p_1$ and $q_1$
 having in $O'(m_0)$ the form respectively 
$$
p_1(m)=\left\{%
\begin{array}{ll}
    0 & \hbox{when  $m<m_0$,} \\
    1 & \hbox{when $m_0<m$,} \\
\end{array}%
\right.\ \
q_1(m)=\left\{%
\begin{array}{ll}
    1 & \hbox{when  $m<m_0$,} \\
    0 & \hbox{when  $m_0<m$,} \\
\end{array}%
\right.
$$
and such that  $p_1(m)+q_1(m)a_{-h}(m)\ne 0$ out of $O'(m_0)$. 

Let us write out the local representatives for this operator:
$$
(P_1+Q_1U_{-h})^{\pm}=\left\{%
\begin{array}{ll}
    a_{-h} & \hbox{ when  $m<m_0,\ m\in O'(m_0)$,} \\
    1 & \hbox{when  $m_0<m,\ m\in O'(m_0)$,} \\
    p_1(m)+q_1(m)a_{-h}(m) & \hbox{when  $m\notin O'(m_0)$,} \\
\end{array}%
\right.
$$
$$
(P_1+Q_1U_{-h})(m_0,+,\pm)=1, \ \
(P_1+Q_1U_{-h})(m_0,-,\pm)=T_{-h},
$$
$$
(P_1+Q_1U_{-h})(m_0)(t)=\begin{pmatrix}
  1 & 0 \\
  0 & 0
\end{pmatrix}+\begin{pmatrix}
  0 & 0 \\
  0 & 1
\end{pmatrix}T_{-h}.
$$

\begin{Tm}
\label{t3} 
 A Fredholm operator  $d$ of the form  (\ref{e2})
satisfies the inequalities \\
(case III)
\begin{equation}
\label{caseIII}
 |b^{\pm}_0(m_0+0)|>|b^{\pm}_1(m_0+0)|,\
    |b^{\pm}_0(m_0-0)|<|b^{\pm}_1(m_0-0)|
\end{equation}
or  the inequalities  (case  IV)
\begin{equation}\label{caseIV}
 |b^{\pm}_0(m_0+0)|<|b^{\pm}_1(m_0+0)|,\
    |b^{\pm}_0(m_0-0)|>|b^{\pm}_1(m_0-0)|)
\end{equation}
iff there exist zero index operators  $w$ 
and  $s$  from  $B$ such that the local representative for 
$e=w^{-1}ds$ on    ${\mathbf {R}}_{m_0}$
is equal to  (\ref{e3}) (is equal to  (\ref{e4})) for  $l=1$. 
Moreover there exists a homotopy in the class of Fredholm operators from 
$C^\ast (B,U_h)$ of the 
operator $e$ to the operator $e'$ such that 
in case III
$$
d' =e'(P_1 +Q_1 U_{-h}) \in B
$$
and in case IV
$$
d' = e' (P_1U_h + Q_1) \in B
$$
and 
$$
\mathrm{ind}\ d =\mathrm{ind}\ d'.
$$
\end{Tm}

{\bf Proof.} The proof of this theorem is  also based on  
theorem \ref{t1} and thus we are using the same notation.
The inequalities 
 (\ref{caseIII}) and  (\ref{caseIV}) follow from the cases (i) and  (ii) of 
 lemma \ref{l2} where one should take the corresponding infinities. 
 The further proof is analogous to that of the preceding theorem. 
 
In  case  III for every  $\tau\in[0,1]$ we set
$$
(e'_{\tau})^{\pm}(m)=\left\{%
\begin{array}{ll}
    (w^{-1})^{\pm}(m)({\varphi}_{\tau}(m)b_0^{\pm}(m) +b_1^{\pm}(m)a_h(m))s^{\pm}(m) & 
    \hbox{when  $m<m_0$,} \\
    (w^{-1})^{\pm}(m)(b_0^{\pm}(m) +{\varphi}_{\tau}(m)b_1^{\pm}(m)a_h(m))s^{\pm}(m) & 
    \hbox{when  $m_0<m$,} \\
\end{array}%
\right.
$$

$$
e'_{\tau}(m_0,+,\pm)=1+(1-\tau)
\frac{b_1^{\pm}(m_0+0)}{b_0^{\pm}(m_0+0)} T_h,
$$
$$
e'_{\tau}(m_0,-,\pm)=(1-\tau)
\frac{b_0^{\pm}(m_0-0)}{b_1^{\pm}(m_0-0)} +T_h,
$$
$$
e'_{\tau}(m_0)(t)=\begin{pmatrix}
  1 & 0 \\
  0 & (1-\tau)e^0_{22}(t)
\end{pmatrix}+\begin{pmatrix}
  (1-\tau)e^1_{11}(t) & 0 \\
  0 & 1
\end{pmatrix}T_h,
$$
where the function  ${\varphi}_{\tau}$ is defined in the same way as in the 
proof of the preceding theorem.

For  $\tau=1$ we have 
$$
(e'_1)^{\pm}(m)=\left\{%
\begin{array}{ll}
    (w^{-1})^{\pm}(m)({\varphi}_1(m)b_0^{\pm}(m) +b_1^{\pm}(m)a_h(m))s^{\pm}(m) & 
    \hbox{when  $m<m_0$,} \\
    (w^{-1})^{\pm}(m)(b_0^{\pm}(m) +{\varphi}_1(m)b_1^{\pm}(m)a_h(m))s^{\pm}(m) & 
    \hbox{when $m_0<m$,} \\
\end{array}%
\right.
$$
$$
e'_1(m_0,+,\pm)=1,\ \ e'_1(m_0,-,\pm)=T_h,
$$
$$
e'_1(m_0)(t)=\begin{pmatrix}
  1 & 0 \\
  0 & 0
\end{pmatrix}+\begin{pmatrix}
  0 & 0 \\
  0 & 1
\end{pmatrix}T_h.
$$

For every  $\tau\in[0,1]$    the invertible operators  $(e'_{\tau})^{\pm}(m)$,
$e'_{\tau}(m_0,+,\pm)$, $e'_{\tau}(m_0,-,\pm)$ and  $e'_{\tau}(m_0)$
define the symbol of a certain Fredholm operator 
$e'_{\tau}$ and 
$$
\mathrm{ind}\ d=\mathrm{ind}\ e=\mathrm{ind}\ e'_0=\mathrm{ind}\
e'_1.
$$
Set $e'=e'_1$. Clearly  the operator  $d'=e'(P_1+Q_1U_{-h})$ belongs to the algebra 
and its index can be calculated by  (\ref{indB}). Thus 
$$
\mathrm{ind}\ d=\mathrm{ind}\ e'=\mathrm{ind}\
(e'(P_1+Q_1U_{-h})).
$$

To reduce the case  IV to the case  III it is enough to multiply the operator 
$d$ by  $U^{-1}_h$ . The proof is complete. $\square $

\

An important technical step in the investigation of the cases 
to be considered is the next result. 

\begin{Tm}
\label{t4}
\it{Let  $P$ and  $Q$ be the projections on  $L^2(M)$ having the form 
$$
P(m)=\frac{1}{2}(I+S),\ \ \ Q(m)=\frac{1}{2}(I-S),
$$
where  $I$ is the identity operator and   $S$ is the singular integral operator.
Then the operators  $P+QU_h$ and  $Q+PU_h$
are Fredholm and their indexes are zero. 
}
\end{Tm}

{\bf Proof.} Since  
$Q+PU_h=(P+QU_{-h})U_h$ it is enough to verify the statement for the operator 
 $P+QU_h$.

An operator is Fredholm iff all its local representatives are invertible. 
The explicit form of the local representatives for  $P+QU_h$ is written out below 
and it is clear that all of them are invertible: 
$$
((P+QU_h)^+)^{-1}(m)=1,\ \ \ ((P+QU_h)^-)^{-1}(m)=a_{-h}(m),
$$
$$
((P+QU_h)(m_0,\pm,+)^{-1}=1,\ \ \
((P+QU_h)(m_0,\pm,-)^{-1}=T_{-h},
$$
$$
((P+QU_h)(m_0)(t))^{-1}=$$ $$=\frac{1}{e^{\pi h}+e^{2\pi t}}\biggl(\left(%
\begin{array}{cc}
  e^{2\pi t} & \mathrm{i}e^{\pi t} \\
  -\mathrm{i}e^{\pi t+\pi h} & e^{\pi h} \\
\end{array}%
\right)+\left(%
\begin{array}{cc}
  e^{\pi h} & -\mathrm{i}e^{\pi t} \\
  \mathrm{i}e^{\pi t+\pi h} & e^{2\pi t} \\
\end{array}%
\right)T_{-h}\biggr).
$$

Let us show now that the index of  $P+QU_h$  is zero. 
Let 
$F_{\alpha}:L^2(M)\to L^2(M)$
be the operator of the form 
$$
\bigl[ F_{\alpha}(f)\bigr](m)=|{\alpha}'(m)|^{1/2}f(\alpha(m)),
$$
where the diffeomorphism  $\alpha:M\to M$ changes the orientation of the curve 
 $M$ and on the neighbourhood  $O(m_0)$ acts as the  reflection 
 and  $\alpha(m_0)=m_0,\
{\alpha}^2(m)=m$. It is known (see [\ref{Litvinchuk}]) that this operator 
possesses  the following properties:
\begin{enumerate}
    \item[(1)] $F^2_{\alpha}=I$,
    \item[(2)] $F_{\alpha}PF_{\alpha}\sim Q$,
    \item[(3)] $F_{\alpha}QF_{\alpha}\sim P$,
    \item[(4)] $F_{\alpha}U_hF_{\alpha}=U_h$
\end{enumerate}
(here the sign  $\sim$ means the the left hand part and the right hand part differ 
by a compact summand).

It follows from (2)---(4) that 
$$
F_{\alpha}(P+QU_h)F_{\alpha}=F_{\alpha}PF_{\alpha}+F_{\alpha}QF_{\alpha}F_{\alpha}U_hF_{\alpha}\sim
Q+PU_h.
$$
Thus 
$$
\mathrm{ind}\ (P+QU_h)=\mathrm{ind}\ (Q+PU_h)=\mathrm{ind}\
(PU^{\ast}_h+Q)= \mathrm{ind}\ (QU^{\ast}_h+P),
$$
and therefore 
$$
2\mathrm{ind}\ (P+QU_h)=\mathrm{ind}\
\bigl[(P+QU_h)(P+QU^{\ast}_h)\bigr]=\mathrm{ind}\ \bigl[
P+QU_hP+QU_hQU^{\ast}_h\bigr],
$$
where $\tau\in[0,1]$.

Consider the operator 
$$
X+(1-\tau)YU_h=\bigl[ P+QU_hQU^{\ast}_h\bigr]+(1-\tau) QU_hP.
$$
Since 
$$
(X+(1-\tau)YU_h)^{\pm}(m)=1,
$$
$$
(X+(1-\tau)YU_h)(m_0,+,\pm)=(X+(1-\tau)YU_h)(m_0,-,\pm)=1,
$$
$$
\mathrm{det}\ (P+QU_hQU^{\ast}_h)(m_0)(t)=\frac{(1+e^{\pi h+2\pi
t})^2}{(1+e^{2\pi t})(1+e^{2\pi h+2\pi t})}\ne0 \ \ \forall
t\in{\mathbf {R}}_{m_0},
$$
$$
X(m_0)(\pm\infty)=1>0=Y(m_0)(\pm\infty).
$$
it follows from theorem \ref{t2} that the operator $X+(1-\tau)YU_h$
is Fredholm.

Thus 
$$
2\mathrm{ind}\ (P+QU_h)= $$ $$=\mathrm{ind}\ \bigl[ P+(1-\tau)
QU_hP+QU_hQU^{\ast}_h\bigr]=\mathrm{ind}\ \bigl[ P+0\cdot
QU_hP+QU_hQU^{\ast}_h\bigr]=$$ $$=-\frac{1}{2\pi}\left( \arg
1|_{M_+}- \arg 1|_{M_-}+\arg\left[\frac{(1+e^{\pi h+2\pi
t})^2}{(1+e^{2\pi t})(1+e^{2\pi h+2\pi t})}\right]\bigg |
^{+\infty}_{t=-\infty} \right)=0.
$$
The theorem is proved. $\square $

\begin{Rk} 
\label{r4} 
{\em  The operator  $(P+QU_h)(m_0)(t)$
satisfies the equality
$$
\widetilde w_1^{-1}(t)\bigg[(P+QU_h)(m_0)(t)\bigg]\widetilde
s_1(t)=e_0(t)+e_1(t)T_h
$$
where 
$$
e_0(t)=\left(%
\begin{array}{cc}
  0 & 0 \\
  0 & 1 \\
\end{array}%
\right),\ \ \
e_1(t)=\left(%
\begin{array}{cc}
  1 & 0 \\
  0 & 0 \\
\end{array}%
\right),
$$
$$
\widetilde w_1(t)=\left(%
\begin{array}{cc}
  1 & - \mathrm{i}e^{\pi t} \\
  \mathrm{i}e^{\pi t} & 1 \\
\end{array}%
\right),
$$
$$
\widetilde s_1(t)=\left(%
\begin{array}{cc}
  \frac{e^{-\pi h}(1+e^{\pi t})(e^{2\pi h}+e^{2\pi t})}{(1+e^{\pi t})(e^{\pi h}+e^{2\pi t})(1+e^{\pi t-\pi h})} &
  \frac{- \mathrm{i}e^{-\pi h+\pi t}(e^{\pi h}+e^{\pi t})(1+e^{2\pi t})}{(1+e^{\pi t})(e^{\pi h}+e^{2\pi t})(1+e^{\pi t-\pi h})}  \\
  \frac{\mathrm{i}e^{-\pi h+\pi t}(e^{2\pi h}+e^{2\pi t})(1+e^{\pi t})}{(1+e^{\pi t})(e^{\pi h}+e^{2\pi t})(1+e^{\pi t-\pi h})} &
  \frac{-(e^{\pi h}+e^{\pi t})(1+e^{2\pi t})}{(1+e^{\pi t})(e^{\pi h}+e^{2\pi t})(1+e^{\pi t-\pi h})}  \\
\end{array}%
\right).
$$
}
\end{Rk}

\begin{Tm}
\label{t5}
 \it{A Fredholm operator  $d$ of the form  (\ref{e2})
satisfies the inequalities  (case  V)
\begin{equation}
\label{caseV}
 |b^{+}_0(m_0\pm 0)|>|b^{+}_1(m_0\pm 0)|,\
    |b^{-}_0(m_0\pm 0)|<|b^{-}_1(m_0\pm 0)|
\end{equation}
or the inequalities  (case VI)
\begin{equation}
\label{caseVI}
 |b^{-}_0(m_0\pm 0)|>|b^{-}_1(m_0\pm 0)|,\
    |b^{+}_0(m_0\pm 0)|<|b^{+}_1(m_0\pm 0)|
\end{equation}
iff there are exist non degenerate continuous matrix functions of order  2 \ $w_1$ and  $s_1$ on 
 $\overline {\mathbf {R}}_{m_0}$ such that the operator  $e=w^{-1}_1d(m_0)s_1$ 
 has the form 
(\ref{e3}) in  case  V or the form  (\ref{e4}) in  case  VI for  $l=1$.
In the case V the matrices   $w(m_0)(+\infty)$ and  $s(m_0)(+\infty)$ are diagonal 
and the matrices  $w(m_0)(-\infty)$ and $s(m_0)(-\infty)$ are skew diagonal while in the case 
 VI we have the opposite situation.

Moreover there exists a homotopy in the class of Fredholm operators of the algebra 
 $C^{\ast}(B, U_h)$ between the operator  $d$ and an operator  $d_1$ and operators of zero indexes 
  $w_2, s_2\in B$ such that in case  V
$$
d'=w_2d_1s_2(PU_h+Q)^{-1}\in B,
$$
and in case  VI
$$
d'=w_2d_1s_2(P+QU_h)^{-1}\in B.
$$
And 
$$
\mathrm{ind\ }d=\mathrm{ind\ }d'.
$$}
\end{Tm}
{\bf Proof.} The proof of this theorem is also based on theorem \ref{t1}. 
The inequalities 
 (\ref{caseV}) and  (\ref{caseVI})
can be obtained from the cases (i) and  (ii) of lemma \ref{l2} by 
taking the corresponding infinities.    

In case  VI we have 
$$
e(t)=w^{-1}_1(t)d(m_0)(t)s_1(t)=\left(%
\begin{array}{cc}
  e^0_{11}(t) & 0 \\
  0 & 1 \\
\end{array}%
\right)+\left(%
\begin{array}{cc}
  1 & 0 \\
  0 & e^1_{22}(t) \\
\end{array}%
\right)T_h,
$$
where  $t\in \overline{\mathbf {R}}_{\mathbf{m_0}}$.


Let us define the homotopy 
$$
e'_{\tau}(t)=\left(%
\begin{array}{cc}
  (1-\tau)e^0_{11}(t) & 0 \\
  0 & 1 \\
\end{array}%
\right)+\left(%
\begin{array}{cc}
  1 & 0 \\
  0 & (1-\tau)e^1_{22}(t) \\
\end{array}%
\right)T_h,\ \tau\in [0,1].
$$
This homotopy defines the homotopy of the operator  $d(m_0)(t)$:
$$
d_{\tau}(m_0)(t)=w_1(t)e'_{\tau}(t)s^{-1}_1(t),
$$
where we also have 
$$
e'_0(t)=e(t),\ d_0(m_0)(t)=d(m_0)(t),
$$
$$
e'_1(t)=\left(%
\begin{array}{cc}
  0 & 0 \\
  0 & 1 \\
\end{array}%
\right)+\left(%
\begin{array}{cc}
  1 & 0 \\
  0 & 0 \\
\end{array}%
\right)T_h,
$$
$$
d_1(m_0)(t)=w_1(t)\bigg[\left(%
\begin{array}{cc}
  0 & 0 \\
  0 & 1 \\
\end{array}%
\right)+\left(%
\begin{array}{cc}
  1 & 0 \\
  0 & 0 \\
\end{array}%
\right)T_h \bigg] s^{-1}_1(t).
$$

The extension of the homotopy to the points  $(m_0,\pm,+)$ and   $(m_0,\pm,-)$
is defined in the following way
$$
d_{\tau}(m_0,\pm,+)=(1-\tau)b_0(m_0,\pm,+)+b_1(m_0,\pm,+)T_h,
$$
$$
d_{\tau}(m_0,\pm,-)=b_0(m_0,\pm,-)+(1-\tau)b_1(m_0,\pm,-)T_h.
$$
Therefore 
$$
d_0(m_0,\pm,+)=d(m_0,\pm,+)=b^+_0(m_0\pm 0)+b^+_1(m_0\pm 0)T_h,
$$
$$
d_0(m_0,\pm,-)=d(m_0,\pm,-)=b^-_0(m_0\pm 0)+b^-_1(m_0\pm 0)T_h,
$$
$$
d_1(m_0,\pm,+)=b^+_1(m_0\pm 0)T_h, \ \ \
d_1(m_0,\pm,-)=b^-_0(m_0\pm 0).
$$

On  $M_+$ and $M_-$  we define the homotopy on the following way
$$
d^+_{\tau}(m)={\varphi}_{\tau}(m)b^+_0(m)+b^+_1(m)a_h(m),
$$
$$
d^-_{\tau}(m)=b^-_0(m)+{\varphi}_{\tau}(m)b^-_1(m)a_h(m),
$$
where the function  ${\varphi}_{\tau}$ is defined in the same way as in the proof of theorem \ref{t2}.

The operators  $(d'_{\tau})^{\pm}(m)$, $d'_{\tau}(m_0,+,\pm)$,
$d'_{\tau}(m_0,-,\pm)$ and $d'_{\tau}(m_0)$ are invertible and for each 
$\tau\in[0,1]$ they define the symbol of a certain Fredholm operator  $d'_{\tau}$ 
and
$$
\mathrm{ind}\ d=\mathrm{ind}\ d_0=\mathrm{ind}\ d_1.
$$

For $t=\pm\infty$ we have the diagonal matrices  $w_2(t)=\widetilde
w_1(t)w^{-1}_1(t)$ and $s_2(t)=$ $=s_1(t)\widetilde s^{-1}_1(t)$,
(here the matrices  $\widetilde w_1(t)$ and  $\widetilde s_1(t)$ are that defined in remark \ref{r4}). 
These matrices can be extended by means of non vanishing continuous functions on 
$M_{\pm}$ in such a way that the indexes of the corresponding operators  $w_2$ and 
$s_2$ will be equal to zero  (see the proof of lemma \ref{l1}).

The operator  $d'=w_2d_1s_2(P+QU_h)^{-1}$ belongs to the algebra  $B$  and 
$$
\mathrm{ind}\ d=\mathrm{ind}\ d'.
$$
Now the index can be calculated by means of  formula  (\ref{indB}).

To reduce  case  V to case  VI it is enough to multiply the operator $d$
by  $U^{-1}_h$. The proof is complete. $\square $

\

\begin{center}
{\bf References }
\end{center}
\renewcommand{\labelenumi}{\bf \theenumi .}
\begin{enumerate}

\item Simonenko I.B. {\em A new general method of investigation of linear 
operator integral equations. I.} \ Izv. AN SSSR, ser. mat. 1965. V. 29, No 3. p. 567-586.
(Russian)  \label{Simonenko1}

\item Simonenko I.B. {\em A new general method of investigation of linear 
operator integral equations. II.}\  Izv. AN SSSR, ser. mat. 1965. V. 29, No 4.
p. 775-782.
(Russian)  \label{Simonenko2}

\item Simonenko I.B., Chin' Ngok Min' {\em  The local method 
in the theory of one dimensional singular integral equations with piecewise 
continuous coefficients. Noetherity.}\  Izdatel'stvo Rostovskogo Universiteta, 1986.
(Russian)
\label{SimonenkoChin}

\item Gohberg I.Ts., Krupnik N.Ja. {\em Introduction into the theory 
of one dimensional singular integral operators.}\  Kishinev,"Shtiintsa", 1973. (Russian) 
 \label{GohbergKrupnik}

\item Akulich E.V., Lebedev A.V. {\em The symbolic calculus 
for singular integral operators with discontinuous oscillating coefficients.}\ 
Dokl. AN Belarusi. 2003. V. 47, No 1. p. 10-14.
(Russian)
\label{N1}

\item  Plamenevsky B.A. {\em Algebras of pseudodifferential operators.} \
M.,1986.
(Russian)
\label{Plamenevski}

\item
 Antonevich A., Lebedev A. {\em Functional differential equations: I. $C^\ast -$theory.} \
Longman Scientific $\&$ Technical, 70, 1994.\label{IPart}

\item
 Antonevich A., Belousov M., Lebedev A. {\em Functional differential
equations: II. $C^\ast -$applications. Part 2.} \ Longman Scientific $\&$
Technical, 95,  1998. 
\label{IIPart2}

\item Gohberg I.Ts., Krupnik N.Ja. {\em On the algebra generated by singular integral operators 
with piecewise continuous coefficients.} \ Funkts. analiz i egho prilozh.
1970. V. 4, No 3. p.
27-38. (Russian)
\label{GohbKrup1}

\item Gohberg I.Ts., Krupnik N.Ja. {\em Singular integral operators with 
piecewise continuous coefficients and their symbols.}\  Izv. AN SSSR, ser. mat.1971.
V. 35,
 No 4. p. 940-964.
 (Russian)
 \label{GohbKrup2}

\item
 Antonevich A., Belousov M., Lebedev A. {\em Functional differential
equations: II. $C^\ast -$applications. Part 1.}\  Longman Scientific $\&$
Technical,  94, 1998. 
\label{IIPart1}

\item Litvinchuk G.S. {\em Boundary value problems and singular integral 
equations with a shift.}\  M. 1977. (Russian)
\label{Litvinchuk}

\end{enumerate}
\end{document}